\pdfoutput=1 
  
\documentclass[a4paper]{article}
\usepackage[latin2]{inputenc}
\usepackage{color}
\usepackage{amssymb}
\usepackage{amsmath}
\usepackage{amscd}
\usepackage{enumerate}
\usepackage{array}
\usepackage{tabularx}
\usepackage{t1enc}
\usepackage[mathscr]{eucal}
\usepackage{graphicx} 
\usepackage{graphics}
\usepackage{paralist}
\usepackage[all]{xy}
\usepackage{url}
\usepackage{multirow} 
\usepackage[authoryear,sort,round]{natbib} 
                                             
\begin{document}
\newtheorem{defi}{Definition}
\newtheorem{tetel}{Theorem}
\newtheorem{felt}{Assumption}

\newcounter{feltet}
\setcounter{feltet}{0} 
\newcommand{\bit}{ \begin{compactitem}}
\newcommand{\eit}{ \end{compactitem} }
\newcommand{\ben}{ \begin{compactenum}[a.)]}
\newcommand{\een}{ \end{compactenum} }
\renewcommand{\labelitemii}{$\cdot$}
\setdefaultleftmargin{1cm}{0.5cm}{0.5cm}{0.5cm}{0.5cm}{0.5cm}
\newcommand{\hely}{\vspace{2mm}}
\newcommand{\kov}{\quad $\Rightarrow$ \quad}
\newcommand{\norm}[1]{\lVert#1\rVert} 
\newcommand{\ass}[1]{ \vspace{2mm} \textbf{Assumption \stepcounter{feltet}#1\arabic{feltet}}}

\title{\bf Weighted bootstrap in GARCH models}
\maketitle
{\bf László Varga }\\
{ Department of Probability Theory and Statistics, Eötvös 
Loránd University, Budapest, Hungary}\\
{\bf E-mail:} { vargal4@math.elte.hu}\medskip
     
{\bf  András Zempléni}\\
{Department of Probability Theory and Statistics, Eötvös 
Loránd University, Budapest, Hungary} \\
{\bf E-mail:} { zempleni@math.elte.hu}\medskip

\begin{abstract}
GARCH models are useful tools in the investigation of phenomena, where volatility changes are prominent features, like most financial data. The parameter estimation via quasi maximum likelihood (QMLE) and its properties are by now well understood. However, there is a gap between practical applications and the theory, as in reality there are usually not enough observations for the limit results to be valid approximations. We try to fill this gap by this paper, where the properties of a recent bootstrap methodology in the context of GARCH modeling are revealed. The results are promising as it turns out that this remarkably simple method has essentially the same limit distribution, as the original estimatorwith the advantage of easy confidence interval construction, as it is demonstrated in the paper.
     
The finite-sample properties of the suggested estimators are investigated through a simulation study, which ensures that the results are practically applicable for sample sizes as low as a thousand. On the other hand, the results are not 100\% accurate until sample size reaches 100 thousands - but it is shown that this property is not a feature of our bootstrap procedure only, as it is shared by the
original QMLE, too.   
\\[0.5cm]
\end{abstract}

\maketitle
           
\nocite{ci2008, bh2004, hkt2004, cgba2011, prr2006, ling2007, luger2011, fz2004, bhatt2012, koja2011}
      
\section{Introduction}

We investigate bootstrap estimation of the parameters of GARCH processes, which are 
known to be able to capture the main stylized facts
of observed financial series. In these models, the conditional variance is 
expressed as a linear function of the squared
past values of the series. 
\begin{defi}  $(X_t)_{t \in \mathbb{Z}}$ is called a GARCH(p,q) process if 
\label{garchdef}
 \begin{flalign}
  \label{garch1} X_t & = \sqrt{h_t} \eta _t   \\
  \label{garch2} h_t & = \omega_0 + \sum\limits _{i=1}^q \alpha_{0i}X^2_{t-i} + \sum\limits _{j=1}^p \beta_{0j}h_{t-j} 
 \end{flalign}
 where $\eta _t$ $(t \in \mathbb{Z})$ are 
i.i.d. \text{(0,1)} random variables, $\omega_0>0, \alpha_{0i} \geq 0,\beta_{0j} \geq 0$ 
 for $i=1,...,q$ and for $j=1,...,p$.
\end{defi}
It defines a stationary process for a well characterized parameter space, 
its most important features are presented in Section 2. 

The most important question in modeling is the parameter estimation. In case of GARCH models,
the QMLE estimation is the most popular one. This assumes Gaussian distribution 
for the observations, providing reasonable approximations even in the case of
other distributions for the innovation $\eta_t$. We conclude Section 2 with presenting
the properties of this estimator. \\
Of course, there are other estimation methods considered in the literature. The oldest and numerically simplest estimation method for GARCH models is the ordinary least squares (OLS). It performs poorer than the QML method and even for ARCH models the method requires moments of order 8 for the original process (\citet{fzbook}, Chapter 6). An other well known method is the least absolute deviations (LAD) estimation, which outperforms the QML estimator if the innovations are Student's $t$ distributed with 3 or 4 degrees of freedom (\citet{py2003}).
\citet{ling2007} proposed a self-weighted QML estimator for the parameters which is close in some aspects to our considerations. There are also several extensions of these estimators, see \citet{bh2004} and \citet{fzbook}, Chapter 9. \\
Section 3 deals with the main objective of this paper, namely the investigation of bootstrap methods. 
Although there are different approaches
for bootstrapping the GARCH models, (these will be explained in more detail in Section 3) we suggest the 
multiplier bootstrap approach recently proposed by \citet{koja2011} for goodness of fit tests for copulas.
This is a simple generalization of the standard bootstrap procedure, where the 
bootstrap sample is denoted by $(\tau_{ni}X_i)$. 
This method is usually called weighted bootstrap and was investigated as early as in the 1990s \citep{bb1995, pw1993}. 

The bootstrap weights $\tau_{ni}$ $(1\leq i \leq n, \ n \geq 1)$  are supposed
to be independent from the process. We show the asymptotic normality of the
bootstrap estimators under conditions, which are fulfilled in the majority of
practical examples. The weighted bootstrap OLS and LAD estimators for AR(1) and ARCH processes were investigated by \citet{bhatt2012}. 
 
Other bootstrap methods for GARCH models proposed in the literature are the residual bootstrap (for instance, see \citet{hy2003}) and the block bootstrap (\citet{ci2008}). These are tools for constructing confidence intervals for the parameters or for functionals of the parameters (\citet{cgba2011}, \citet{luger2011}, \citet{prr2006}) and for evolving goodness-of-fit tests (\citet{luger2011}, \citet{hkt2004}). Bootstrap methods are especially needed if the errors are heavy-tailed and this is the case in most financial applications. 

Section 4 presents the results of a simulation study, where for simplicity we focus on ARCH(1) models. Here we also investigate the small-sample properties of the QMLE estimators, together with their bootstrap counterparts.
This approach is practical as both the similarities and differences can be demonstrated.  
We give the conclusions in Section 5. The proofs can be found in the Appendix.

\section{GARCH models}

In this Section we summarize the needed fundamentals from the theory of GARCH processes (see \citet{fzbook} for example).

We denote the parameter vector by 
\begin{equation*}
\theta=(\theta_1,...,\theta_{p+q+1})^T=(\omega,\alpha_{1},...,\alpha_{q},\beta_{1},...,\beta_{p})^T,
\end{equation*}
which belongs to the parameter space $\Theta = (0,\infty) \times [0,\infty)^{p+q}$. \\
The true value of the parameters, 
$\theta_0=(\omega_0,\alpha_{01},...,\alpha_{0q},
\beta_{01},...,\beta_{0p})^T$ is unknown. 

\begin{tetel}
If there exists a GARCH($p$,$q$) process (\ref{garch1}) - (\ref{garch2}), 
which is second-order stationary, and if $\omega>0$, then
\begin{equation} \label{stac}
\sum\limits_{i=1}^q \alpha_i + \sum\limits_{j=1}^p \beta_j <1 .
\end{equation}
If (\ref{stac}) holds, the unique strictly stationary solution of model (\ref{garch1}) - (\ref{garch2}) is a weak white
noise.
\end{tetel}

\begin{defi}
Let $(B_t)_{t \in \mathbb{Z}}$ be a strictly stationary sequence of random matrices, and 
E$\left( \log^+ \norm{B_t} \right)< \infty$. The (top) Ljapunov exponent of the sequence $(B_t)_{t \in \mathbb{Z}}$ is 
\begin{equation*}
\lambda := \underset{t \to \infty}{\lim} \frac{1}{t} E \left( \log \norm{B_tB_{t-1} \dots B_1} \right).
\end{equation*}
\end{defi}

The GARCH($p$,$q$) process can be written in vector representation 
\begin{equation*}
\underline{z}_t=\underline{b}_t+A_t  \underline{z}_{t-1},
\end{equation*}
where 

\ $A_t= \left(
\begin{matrix}
\alpha_1\eta_t^2 & & \cdots & & \alpha_q  \eta_t^2 & \beta_1 \eta_t^2 & & \cdots & & \beta_p \eta_t^2 \\
1 & 0 & \cdots & & 0 & 0 & & \cdots && 0  \\
0 & 1 & \cdots & & 0 & 0 & & \cdots && 0  \\
\vdots & \ddots & \ddots && \vdots & \vdots & \ddots & \ddots && \vdots \\
 0 &  & \cdots & 1 & 0 & 0 & & \cdots & & 0  \\
\alpha_1 & & \cdots & & \alpha_q & \beta_1 & & \cdots & & \beta_p \\
0 && \cdots && 0 & 1 & 0 &\cdots && 0 \\
0 && \cdots && 0 & 0 & 1 &\cdots && 0 \\
\vdots & \ddots & \ddots && \vdots & \vdots & \ddots &\ddots && \vdots \\
0 && \cdots && 0 & 0 & 0 &\cdots & 1 & 0 \\
\end{matrix}
\right)\in \mathbb{R}^{(q+p) \times (q+p)}$, 
\hely 
\ $\underline{b}_t= \left( 
\begin{array}{c} 
 \omega \eta_t^2 \\ 0 \\ \vdots \\ 0 \\ \omega \\ 0 \\ \vdots \\ 0 
\end{array}
 \right) \in \mathbb{R}^{q+p}$, \quad
$\underline{z}_t= \left(
\begin{array}{c} 
 X_t^2 \\ X_{t-1}^2 \\ \vdots \\ X_{t-q+1}^2 \\
 \sigma_t^2 \\ \sigma_{t-1}^2 \\ \vdots \\ \sigma_{t-p+1}^2 \end{array}
 \right) \in \mathbb{R}^{q+p}$.

\begin{tetel}
Let $\lambda$ denote the Ljapunov exponent of the matrix sequence $(A_t)_{t \in \mathbb{Z}}$. Then
\begin{equation*}
\lambda <0 \Longleftrightarrow \text{there exists a strictly stationary solution of the GARCH}(p,q)\text{ model}. 
\end{equation*}
\end{tetel}

The following theorem shows that the Ljapunov exponent -- thus the strict stationarity -- is in connection with the existence of moments of the GARCH process, which will be helpful to verify the main results.
\begin{tetel}
\label{tethat}
Let $\lambda$ denote the Ljapunov exponent of the matrix sequence $(A_t)_{t \in \mathbb{Z}}$. Then 
\begin{equation*}
\lambda <0 \Longrightarrow \exists s>0, E \sigma^{2s}< \infty, EX_t^{2s}< \infty
\end{equation*}
where $X_t$ is the strictly stationary solution of the GARCH($p$,$q$) model.
\end{tetel}

From now on we will concentrate on the maximum likelihood estimation of the parameters. Assume that $ \{ x_1,\dots,x_n\} $ are observations from a GARCH(p,q) process 
(strictly stationary solution of the model). 
The Gaussian quasi-likelihood function, conditional on the $x_{1-q},...,x_0,\tilde{\sigma}_{1-p}^2,...,\tilde{\sigma}_0^2$  
initial values, is 
\begin{equation*}
L_n(\theta)= L_n(\theta;x_1,...,x_n)= 
\prod\limits_{i=1}^n \frac{1}{\sqrt{2 \pi \tilde{\sigma}_t^2}} e^{-\frac{x_t^2}{2 \tilde{\sigma}_t^2}}.
\end{equation*}
where the $(\tilde{\sigma}_t^2)_{t \geq 1}$ are recursively defined by the following equation:
\begin{equation*}
\tilde{\sigma}_t^2=\tilde{\sigma}_t^2(\theta)=
 \omega + \sum\limits_{i=1}^q \alpha_i x^2_{t-i} + \sum\limits_{j=1}^p \beta_j \tilde{\sigma}_{t-j}^2(\theta )
\end{equation*}

The QMLE of $\theta$ is defined as the solution $\hat{\theta}_n$ of 
\begin{equation} \label{qmle}
\hat{\theta}_n = \underset{\theta \in \Theta}{\text{argmax}} \ L_n(\theta).
\end{equation}

To maximize the Gaussian likelihood function, we have to minimize the following function:
\begin{equation*}
I_n(\theta )  = \frac{1}{n} \sum\limits _{t=1}^n l_t(\theta ), \quad \text{where} \quad  
l_t(\theta ) = \frac{x_t^2}{\tilde{\sigma} ^2_t (\theta) } +  \log ( \tilde{\sigma}_t^2 (\theta) ).
\end{equation*}

Let $\mathcal{A}_{\theta}(z)$ and $\mathcal{B}_{\theta}(z)$ be the 
generating functions 
\begin{equation*} 
\mathcal{A}_{\theta}(z)= \sum\limits_{i=1}^q \alpha_iz^i,
\end{equation*} 
\begin{equation*} \
\mathcal{B}_{\theta}(z)= 1-\sum\limits_{j=1}^p \beta_jz^j.
\end{equation*} 

The following assumptions \textbf{A1-A6} are sufficient for the quasi-maximum likelihood estimator to have a Normal limit
distribution (see \citet{fz2004}):
\begin{itemize}
\item[\textbf{ A1:}] $\theta_0 \in \Theta \text{ and } \Theta \text{ is compact}$
\item[\textbf{ A2:}] $\gamma(A_0) <0 \text{ and for all } \theta \in \Theta, \sum\limits_{j=1}^p \beta_j <1 $
\item[\textbf{ A3:}] $\eta_t^2 \text{ has a nondegenerate distribution and } E\eta _t^2=1 $
\item[\textbf{ A4:}] $ \text{If } p>0,\mathcal{A}_{\theta_0}(z) \text{ and } \mathcal{B}_{\theta_0}(z) \text{ have no common roots, } \\
 \mathcal{A}_{\theta_0}(1) \neq 0, \alpha_{0q}+\beta_{0p} \neq 0 $
\item[\textbf{ A5:}] $\theta_0 \in \text{int}(\Theta)$
\item[\textbf{ A6:}] $\kappa_{\eta}=E\eta _t^4 < \infty $.
\end{itemize}

\begin{tetel} \label{konziszt}
Let $(\hat{\theta}_n)_{n \geq 1} $ be a sequence of QMLEs satisfying (\ref{qmle}), with initial conditions
\begin{flalign}
\label{kef}
x_{1-q}^2=...=x_0^2= x_1 \qquad 
\tilde{\sigma}_0^2=...=\tilde{\sigma}_{1-p}^2 =x_1^2.
\end{flalign}
Under assumptions \textbf{A1}-\textbf{A4}
\begin{equation*}
\hat{\theta}_n \xrightarrow[n \to \infty]{a.s.} \theta_0.
\end{equation*}
\end{tetel}

\begin{tetel}
\label{hatelo}
Under assumptions \textbf{A1}-\textbf{A6}
\begin{equation}
\label{eq:limdist}
\sqrt{n}(\hat{\theta}_n-\theta_0)  \xrightarrow[n \to \infty]{d} N(0,(\kappa_{\eta}-1)J^{-1}),
\end{equation}
where 
\begin{equation}
\label{Jdef}
J:= E_{\theta_0} \left( \frac{\partial^2 l_t(\theta_0)}{\partial \theta \partial \theta^T} \right) =
E_{\theta_0} \left( \frac{1}{\sigma_t^4(\theta_0)} \frac{\partial \sigma_t^2(\theta_0)}{\partial \theta} 
\frac{\partial \sigma_t^2(\theta_0)}{\partial \theta^T}\right).
\end{equation}
\end{tetel}

With different assumptions, Theorem \ref{konziszt} was first proved by \citet{bhk2003}.
Theorem \ref{hatelo} was proved by \citet{bhk2003} and  by \citet{hy2003}.
\citet{hy2003} also generalized the result to the case in which $E\eta ^4 = \infty$ and the distribution of $\eta^2$ is in the domain of attraction of a Gaussian or stable law with exponent $\zeta \in [1,2)$.

\section{Bootstrap methods}

\subsection{Weighted bootstrap}

We define the bootstrap weights as a triangular sequence of random variables $\tau_{ni}$ $(1\leq i \leq n, \ n \geq 1)$ 
independent from the process:\\
$\begin{array}{ccccc}
\tau_{11}  \\
\tau_{21} & \tau_{22} \\
\vdots & \vdots & \ddots \\
\tau_{n1} &  \tau_{n2} & \dots & \tau_{nn}   \\
\vdots & \vdots &&& \ddots 
\end{array}$

To verify the main results, we need some natural assumptions \textbf{B1-B6} for the bootstrap weights:
\begin{itemize}
\item[\textbf{ B1: }] the weights are independent from the GARCH process
\item[\textbf{ B2: }] $P(\tau _{ni} \geq 0)=1  \quad  1\leq i \leq n, \ n \geq 1 $
\item[\textbf{ B3: }] for all $n$, the first four moments of $\tau_{n1}$, $\dots$, $\tau_{nn}$ are finite and equal
\item[\textbf{ B4: }] $\underset{n \to \infty}{\lim} E \tau _{ni} = 1 \quad i=1,2,... $
\item[\textbf{ B5: }] $\gamma:=\underset{n \to \infty}{\lim} E \tau_{ni}^2 < \infty \quad i=1,2,... $
\item[\textbf{ B6: }] $R(\tau _{ni}^2,\tau _{nj}^2) \xrightarrow[n \to \infty]{}0$ \quad if $i \neq j$. 
\end{itemize}

The usual bootstrap procedure (corresponding to a multinomial distribution) provides a 
suitable choice for weights, as it satisfies the six assumptions above. \\
This holds for the following weights as well (we shall use the first two in the paper): 
\begin{equation*}
(\tau_{n1},...,\tau_{nn}) \sim \text{ Multinom} \left( n;\frac{1}{n},...,\frac{1}{n} \right),
\end{equation*}
\begin{equation*}
(\tau_{n1},...,\tau_{nn}) \sim \text{ i.i.d. Exp(1)},
\end{equation*}
\begin{equation*}
(\tau_{n1},...,\tau_{nn}) \sim \text{ i.i.d. } \Gamma(n,n).
\end{equation*}

Calculating the Gaussian likelihood function for the weighted sample, we get the following modified negative loglikelihood function, to be minimized:
\begin{equation*}
I_n^*(\theta )  = \frac{1}{n} \sum\limits _{t=1}^n l_{nt}^*(\theta ), \quad \text{where} \quad  
l_{nt}^*(\theta ) = \tau_{nt} \left( \frac{x_t^2}{\tilde{\sigma} ^2_t(\theta)} +  \log ( \tilde{\sigma} ^2_t(\theta) ) \right)
\end{equation*}
For example if the weights are $(1,2,0,1,...,1)$ then the second element of the sample is taken twice but the third one is omitted etc.

The bootstrap QMLE of the parameter $\theta$ is defined as the solution $\hat{\theta}_n^*$ of 
\begin{equation} \label{bqmle}
\hat{\theta}_n^* = \underset{\theta \in \Theta}{\text{argmax}} \ I_n^*(\theta).
\end{equation}

\begin{tetel}
\label{tetelboot1}
Let $(\hat{\theta}^*_n)_{n \geq 1} $ be a sequence of bootstrap QMLEs satisfying (\ref{bqmle}), with initial conditions (\ref{kef}). Under assumptions \textbf{A1-A4} and \textbf{B1-B4} 
\begin{equation*}
\hat{\theta}_n^* \xrightarrow[n \to \infty]{a.s.} \theta_0.
\end{equation*}
\end{tetel}

\begin{tetel}
\label{tetelboot2}
Under assumptions \textbf{A1-A6} and \textbf{B1-B6}
\begin{equation}
\label{eq:limdistw}
\sqrt{n}(\hat{\theta}_n^*-\theta_0)  \xrightarrow[n \to \infty]{d} N\left( 0,\gamma (\kappa_{\eta}-1)J^{-1} \right)
\end{equation}
where 
\begin{flalign*}
J:= E_{\theta_0} \left( \frac{\partial^2 l_t(\theta_0)}{\partial \theta \partial \theta^T} \right) =
E_{\theta_0} \left( \frac{1}{\sigma_t^4(\theta_0)} \frac{\partial \sigma_t^2(\theta_0)}{\partial \theta} 
\frac{\partial \sigma_t^2(\theta_0)}{\partial \theta^T}\right).
\end{flalign*}
\end{tetel}

The proofs of Theorems \ref{tetelboot1} and \ref{tetelboot2} can be found in the Appendix.

\subsection{Residual bootstrap}
    
A residual bootstrap method was proposed by \citet{hy2003}, who also constructed one-sided bootstrap confidence intervals and analyzed its coverage percentages by simulations for stationary ARCH(2) and GARCH(1,1) processes.\\
The construction of the residual bootstrap sample consists of the following steps, which turns out to be useful if
the sample is in its stationary distribution and we apply a suitable burn-in period:

\begin{enumerate}
 \item Given a sample $\{ x_1,...,x_n \}$, compute the QMLE $\hat{\theta}_n$:
  \begin{center}
  $\hat{\theta}_n=  \underset{\theta \in \Theta}{\text{argmin}} \frac{1}{n} \sum\limits _{t=1}^n l_t(\theta )$.
  \end{center}
 \item Estimate the conditional variance $\hat{\sigma}_t$ of the process
  \begin{center}
  $\hat{\sigma}_t= \sqrt{ \tilde{\sigma}_t^2 (\hat{\theta}_n) }$ \quad $t=1,...,n$. 
  \end{center}
 \item Estimate the residuals $\tilde{\eta}_t$
  \begin{center}
  $\tilde{\eta}_t= \frac{x_t}{\hat{\sigma}_t}$ \quad $t=1,...,n$.
  \end{center}
 \item Calculate the standardized residuals $\hat{\eta}_t$
  \begin{center}
  $\hat{\eta}_t= \frac{\tilde{\eta}_t-\frac{\sum_s \tilde{\eta}_s}{n}}
   {\sqrt{ \frac{\sum_s \tilde{\eta}_s^2}{n} - \left(\frac{\sum_s \tilde{\eta}_s}{n}\right)^2 }}$ \quad $t=1,...,n$. 
  \end{center}
 \item Draw a bootstrap sample with replacement from the standardized residuals: 
   $\{ \eta^*_1,...,\eta_n^* \}$.
 \item Using $\hat{\theta}_n$ and $\{ \eta^*_1,...,\eta_n^* \}$, let us compute the residual boostrap sample 
   $\{ x^*_1,...,x_n^* \}$ of the process
   \begin{flalign*}
   x_t^* & = \sigma^*_t\eta _t^*  \qquad t=1,...,n  \\
   \left(\sigma_t^* \right)^2 & = \hat{\omega} + \sum\limits _{i=1}^q \hat{\alpha}_i(x^*_{t-i})^2 +
     \sum\limits_{j=1}^p \hat{\beta}_j (\sigma^*_{t-j})^2. 
  \end{flalign*}
\end{enumerate}

By means of this residual bootstrap procedure, also confidence intervals for future values of the time series and for the $\sigma_t$ volatilities can be constructed \citep{prr2006}.

\section{Simulations}

\begin{figure}[ht]
  \includegraphics[width=12cm]{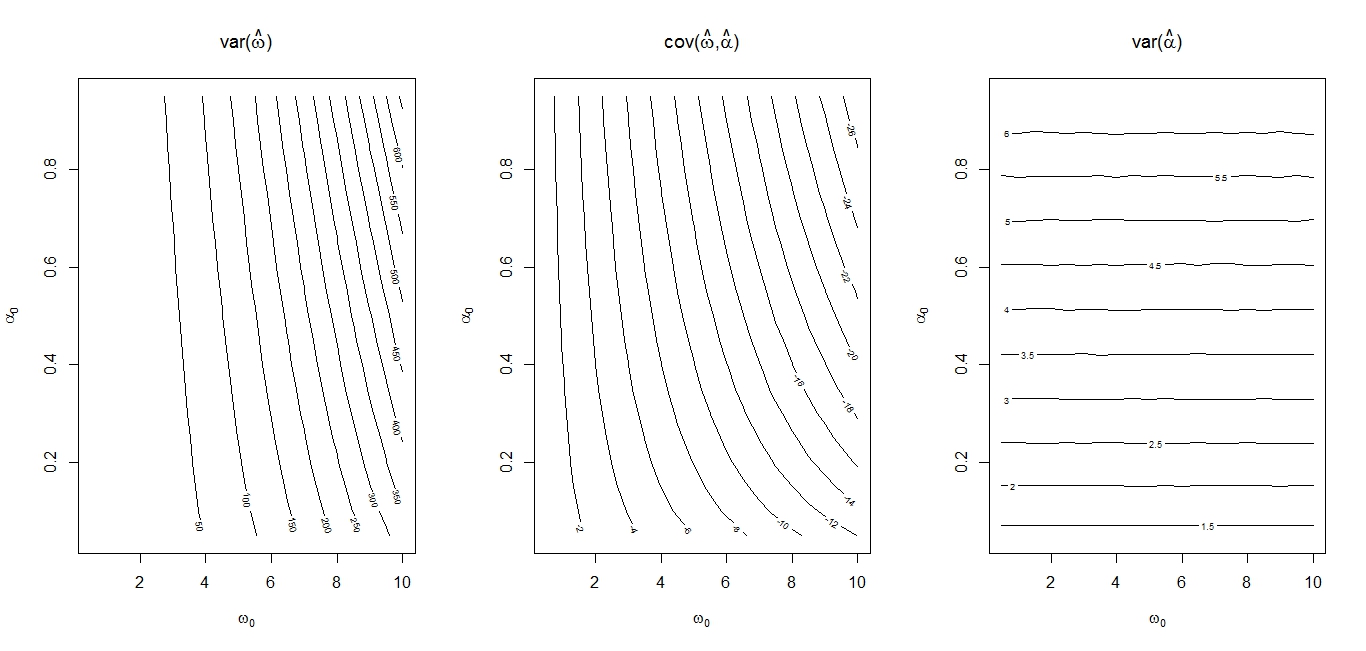} 
  \caption{Contours of the elements of the limiting covariance matrix, ARCH(1) process}
  \label{fig:kontur} 
\end{figure}

Although the GARCH(1,1) models perform usually better and surprisingly well against other, more sophisticated models \citep[see][]{hl2005}, for the sake of simplicity we decided to illustrate the main results for stationary ARCH(1) models (special case $p$=0, $q$=1 of Definition \ref{garchdef}). So suppose that $(X_t)_{t \in \mathbb{Z}}$  is generated by the ARCH(1) process
 \begin{flalign*}
 X_t & = \sqrt{\omega_0 + \alpha_{0}X^2_{t-1}}\eta _t  , 
 \end{flalign*}
 where $\eta _t$ $(t \in \mathbb{Z})$ are 
i.i.d. \text{(0,1)} random variables, and $\theta_0=(\omega_0,\alpha_{0})$, $\omega_0>0$, $\alpha_{0}\geq 0$ are the true parameters. 
The covariance matrix $(\kappa_{\eta}-1)J^{-1}$ of the limit distribution of the QMLE depends on the true parameters. We analyzed this dependence in stationary ARCH(1) processes, where the parameters are $\omega_0>0$ and $0< \alpha_0 <1$. 
The matrix $J$ itself can only be approximated via simulations derived from (\ref{Jdef}): for large $N$ and simulated data $(x_t)_{t=1,...,N}$, 
\begin{equation*}
J \approx \hat{J} = \frac{1}{N}\sum\limits_{t=1}^N \frac{1}{(\omega_0+\alpha_0x^2_t)^2}
 \left( \begin{matrix} 1 \\ x^2_t \end{matrix} \right) \left( \begin{matrix} 1 & x^2_t \end{matrix} \right).
\end{equation*} \\
Figure \ref{fig:kontur} displays the contours of the elements of the limiting covariance matrix if the innovations are Gaussian, based on $N=10^8$ simulations, which provides accurate results up to at least four digits. The variance of $\hat{\omega}$ and the covariance between $\hat{\omega}$ and $\hat{\alpha}$ are both more sensitive to changes in $\omega_0$ than in $\alpha_0$. The variance of the estimated parameter $\hat{\alpha}$ does not seem to depend on the true parameter value $\omega_0$. This is not trivial from the theoretical results, as from (\ref{Jdef}) we get
\begin{flalign*}
\text{var}(\hat{\alpha})=\frac{ E_{\theta_0} \left( \frac{1}{\omega_0+\alpha_0X^2_{t-1}} \right)} 
{E_{\theta_0} \left( \frac{1}{\omega_0+\alpha_0X^2_{t-1}} \right) 
E_{\theta_0} \left( \frac{X^4_{t-1}}{\omega_0+\alpha_0X^2_{t-1}} \right) -
E_{\theta_0}^2 \left( \frac{X^2_{t-1}}{\omega_0+\alpha_0X^2_{t-1}} \right)},
\end{flalign*}
which needs further investigation.
  
\begin{figure}[ht] 
 \begin{center} \includegraphics[width=6cm]{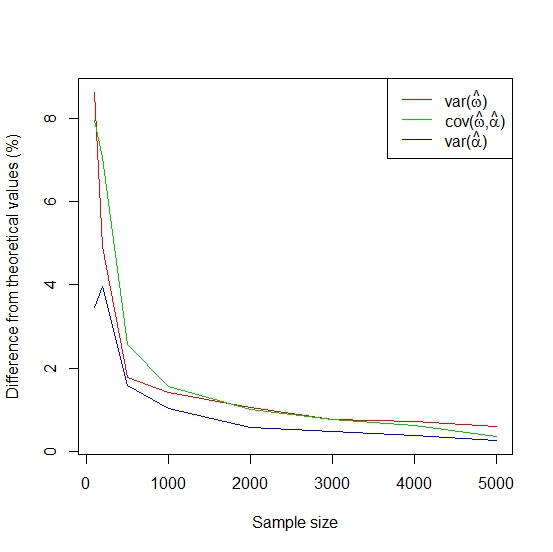}  \end{center} 
 \caption{Convergence of the sample covariance matrix, ARCH(1) process, $\omega_0=1$ and $\alpha_0=0.5$ }
 \label{fig:konv1}
\end{figure}

From now on we will concentrate on the ARCH(1) process with parameters $\omega_0=1$ and $\alpha_0=0.5$.
Then the limiting covariance matrix of the QML estimation is 
\begin{equation} \label{theomat}
\left( \begin{matrix} 4.893 & -2.148	 \\ -2.148	  & 3.926 \end{matrix} \right).
\end{equation}
Unfortunately (minimum) $10^6$ replications are needed to confidently estimate the matrix, which 
takes several hours for an i7 computer with 8 GB RAM memory. We will see that even the bootstrap can't help much if we draw too few samples. \\
We drew $10^6$ samples with Gaussian innovations of size 100 to 5000 and calculated the covariance matrix of the QML estimations. 
Figure \ref{fig:konv1} shows that the rate of convergence drastically improves until the sample size is under 1000 and just slightly after that.
We found also for other pairs of parameters that with simulations of sample size 2000, the covariance matrix can be estimated quite well, within a 1\% margin. 

\begin{figure}[ht] 
 \begin{center}  \includegraphics[width=12cm]{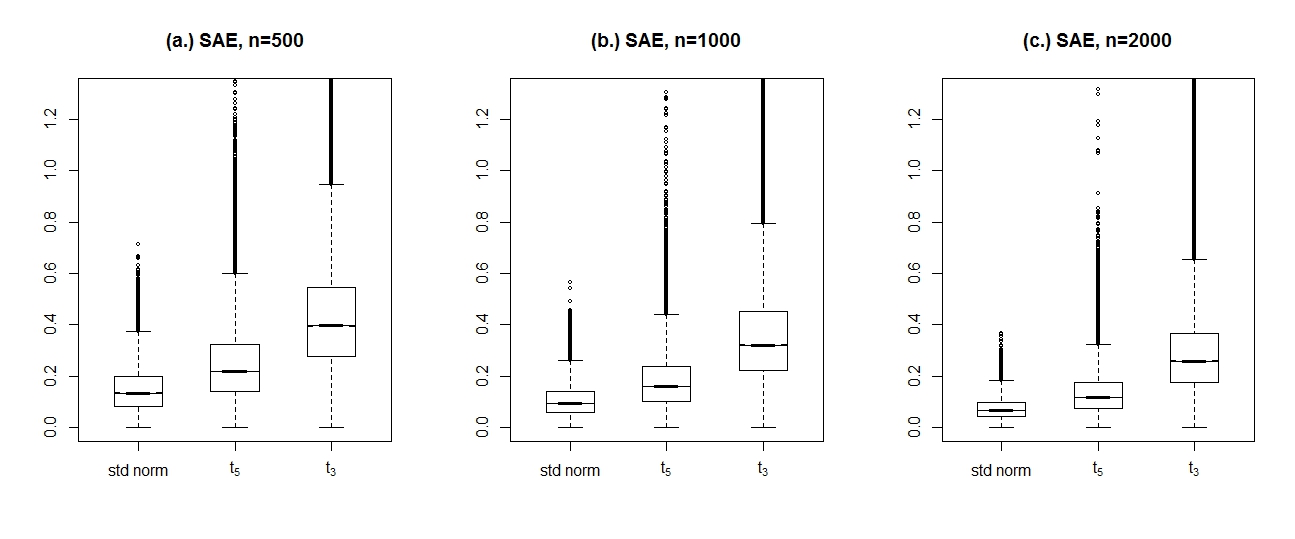} \end{center} 
  \caption{Boxplots of the sum of absolute errors (SAE) of the parameters if the innovations are standard Gaussian, Student's $t$ with 5 and 3 degrees of freedom for different sample sizes: (a.) n=500; (b.) n=1000; (c.) n=2000}
  \label{fig:sae}
\end{figure}

After that, 50000 samples of size n=500, 1000 and 2000 were generated with standard Gaussian and Student's $t$ distributed innovations with 5 and 3 degrees of freedom, and we estimated the parameters with 
the QML method, described in Section 2. Boxplots of the sum of absolute errors (SAE) are depicted in Figure \ref{fig:sae}. The SAE is defined as $|\hat{\omega}-\omega_0|+|\hat{\alpha}-\alpha_0|$. We can see that the heavier tailes the innovations have, the larger the SAE is. Note that the Student's $t$ errors with 3 degrees of freedom have infinite fourth moment -- so Theorem \ref{hatelo} does not work --, but the quasi maximum likelihood estimates are fairly close on average to the original parameters. As the sample size increases, the SAEs of course become smaller. Figure \ref{fig:sae} doesn't display all SAE values for the Student's $t$ innovations, 
the results for some samples are so bad 
that the SAE of the estimated parameters is more than 100.

\begin{figure}[ht] 
 \begin{center}  \includegraphics[width=6cm]{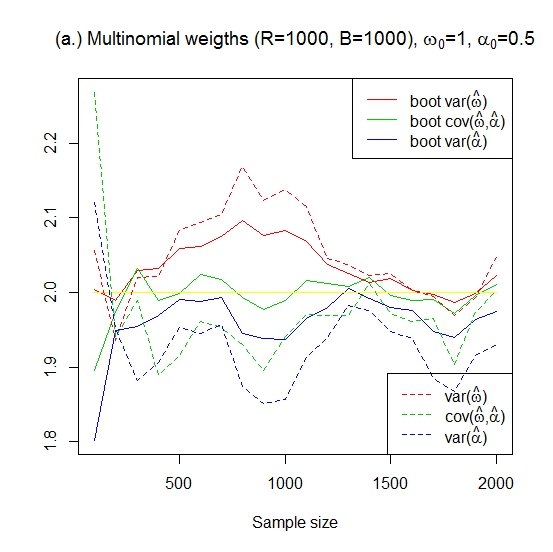}  
 \includegraphics[width=6cm]{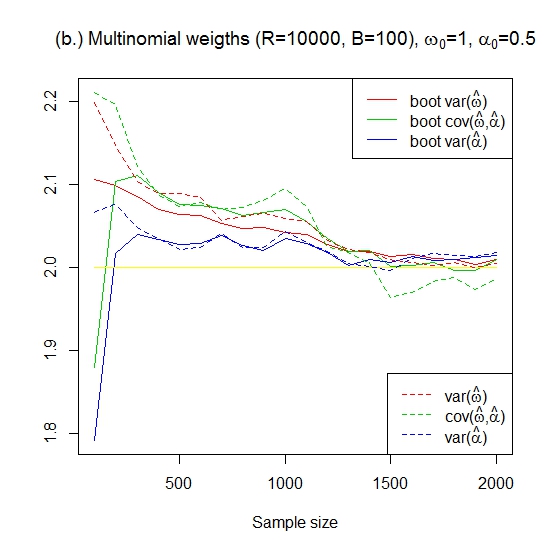} \\
 \includegraphics[width=6cm]{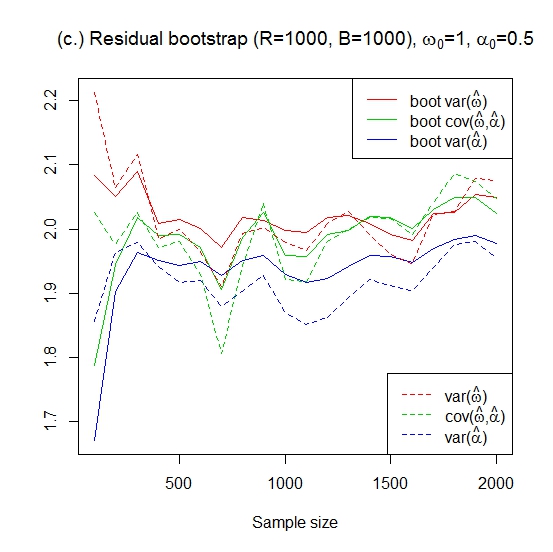}
 \includegraphics[width=6cm]{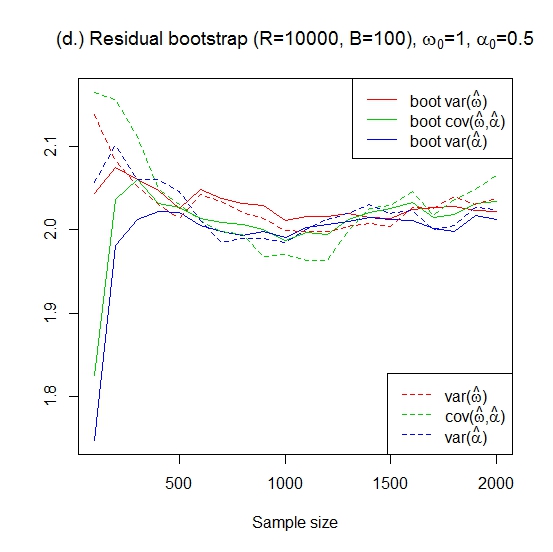} \end{center} 
 \vspace{-5mm}
  \caption{Convergence of the sample covariance matrix, ARCH(1) process, $\omega_0=1$ and $\alpha_0=0.5$; (a.) Weighted bootstrap with multinomial weights, $R$=1000 and $B$=1000; (b.) Weighted bootstrap with multinomial weights, $R$=10000 and $B$=100; (c.) Residual bootstrap, $R$=1000 and $B$=1000; (d.) Residual bootstrap, $R$=10000 and $B$=100}
  \label{fig:konv2}
\end{figure}

If we take multinomially distributed weights, then the scaling factor of the covariance matrix is $\gamma=\underset{n \to \infty}{\lim} E \tau_{ni}^2= \underset{n \to \infty}{\lim} \left( 2- \frac{1}{n} \right) = 2$, therefore the quotient of the two matrices by its elements must be near 2.
Figure \ref{fig:konv2} displays the convergence of the elements of the sample covariance matrix, divided element-wise by the theoretical covariance matrix (\ref{theomat}), if the sample matrices are calculated with the multinomially weighted bootstrap (panels (a.) and (b.)) or with the residual bootstrap (panels (c.) and (d.)), for sample sizes ranging from 100 to 2000. Panel (a.) and (c.) show the convergence based on $R=1000$ samples which were bootstrapped $B=1000$ times, while the other two panels display simulations with $R=10000$ and $B=100$. The dashed lines are the sample covariance matrix values without bootstrap weights, divided by the theoretical values and scaled to 2. \\ 
Unfortunately in Theorem \ref{hatelo} there is not a swift convergence. In panel (a.) of Figure \ref{fig:konv2} we can't see a straight convergence, the bootstrap can't substantially improve the properties of the original samples, it only decreases the differences.  
Panel (b.) of Figure \ref{fig:konv2} helps to understand the reason: the number of samples $R=1000$ was too few. If we raise the number of samples to $R=10000$, and (for practical reasons) 
decrease the bootstrap repetitions to $B=100$, the convergence becomes quite good. 
Looking at the simulations it is not obvious which of the two bootstrap methods is the better one.
     
After that, we constructed 95\% confidence intervals for the GARCH parameters with Gaussian innovations.
Table \ref{tab:tabl} contains the average coverage percentage of confidence intervals for the parameters $\omega$ and $\alpha$ for different sample sizes (500, 1000, 2000) and using residual or weighted bootstrap methods, always compared to the Monte Carlo empirical confidence intervals. 
For sample size of $n=500$, the residual bootstrap outperformed the weighted bootstrap; 
but for sample size 2000, the residual bootstrap performed mostly better then the residual bootstrap. Using the weighted bootstrap, the average coverage of the confidence intervals improved by increasing the sample size which can't be stated in case of residual bootstrap.

\begin{table}[ht]
\begin{center} 
\begin{tabular}{|p{1cm}|c|cc|cc|cc|}
\hline
 Sample size &  \multirow{2}{*}{Method} & \multicolumn{2}{p{2.2cm}|}{Average coverage} & 
  \multicolumn{2}{p{2.2cm}|}{Average coverage below} & \multicolumn{2}{p{2.2cm}|}{Average coverage above} \\ 
 && $\omega$ & $\alpha$ &  $\omega$ & $\alpha$ & $\omega$ & $\alpha$ \\ \hline 
& Monte Carlo & 95\% & 95\% & 2,5\% & 2,5\% & 2,5\% & 2,5\% \\ \hline
\multirow{2}{*}{500} & RB & 94.93 & 95.07 & 2.12 & 2.35 & 2.94 & 2.58 \\   
 & WB & 94.19 & 94.23 &  2.66 & 2.29 & 3.15 & 3.47 \\ \hline
\multirow{2}{*}{1000} & RB & 95.47 & 95.52 &  2.61 & 1.99 & 1.92 & 2.49 \\
 & WB & 94.81 & 94.88 & 3.06 & 1.93 &  2.13 & 3.19 \\ \hline
\multirow{2}{*}{2000} & RB & 95.29 & 94.77 & 2.26 & 2.18 & 2.44 & 3.06 \\
 & WB &  94.74 & 95.07 & 2.62 & 2.21 & 2.64 & 2.72 \\ \hline
\end{tabular} \end{center} 
\caption{Average coverage 
percentages of confidence intervals for the parameters $\omega$ and $\alpha$ for sample sizes 500, 1000, 2000 and using residual bootstrap (RB) or weighted bootstrap (WB) methods.} 
\label{tab:tabl}
\end{table}  

\begin{figure}[ht] 
 \begin{center}   
 \includegraphics[width=12cm]{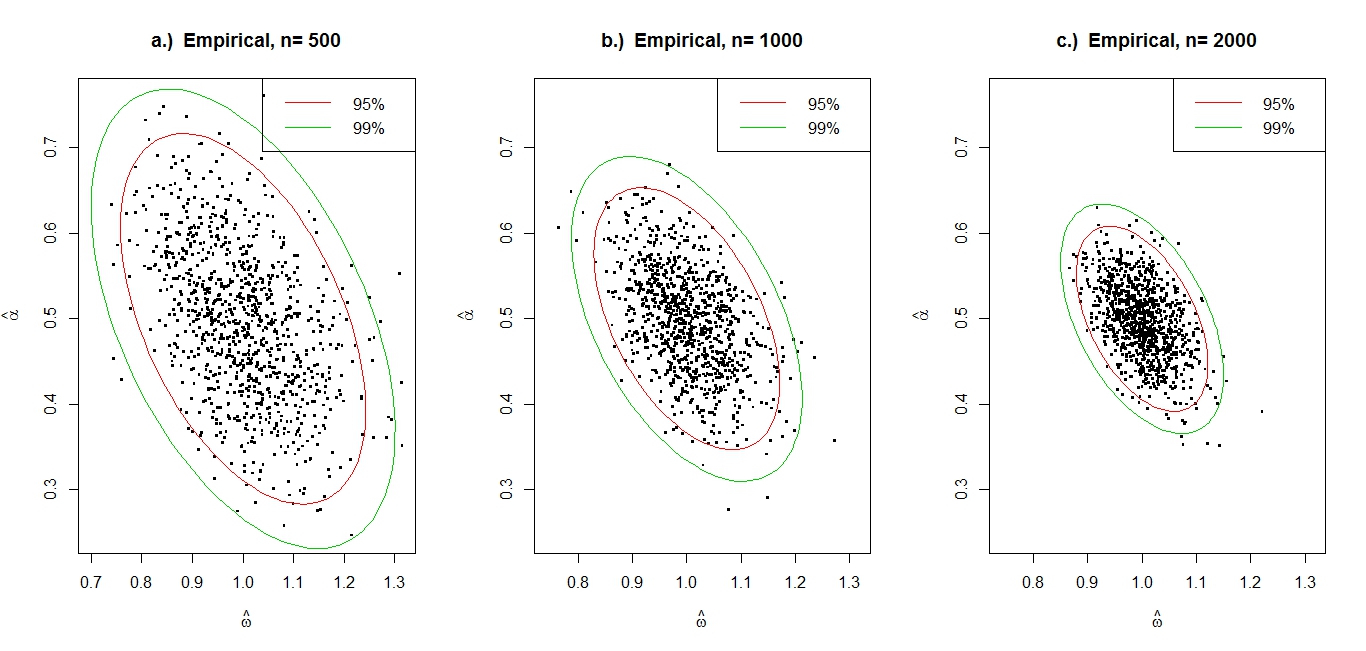} \\
 \includegraphics[width=12cm]{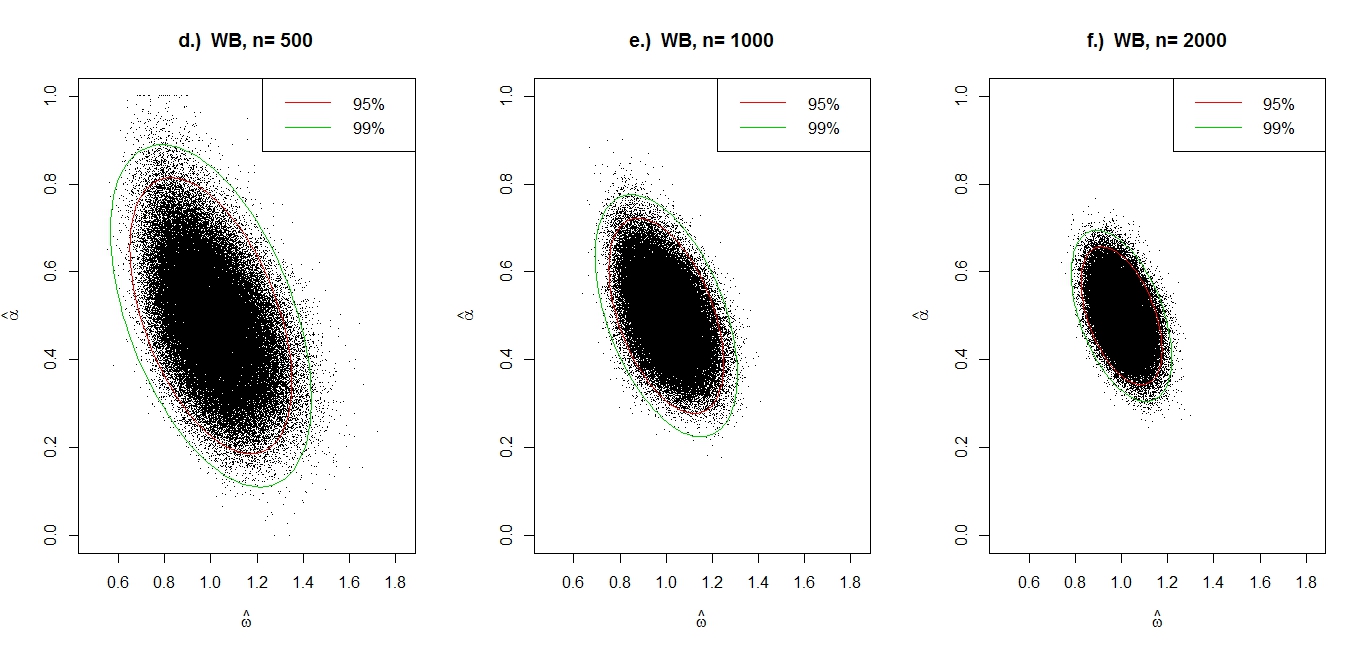} 
 \end{center} 
 \caption{Pairs of estimated parameters ($\hat{\omega},\hat{\alpha}$) and the 95\% and 99\% confidence sets -- according to the limiting distribution -- for different sample sizes; $R=1000$ empirical estimated parameters (a.), b.), c.)) and weighted bootstrap (WB) estimators (d.), e.), f.), $B=100$).}
 \label{fig:confset}
\end{figure}

Using the limiting distributions (\ref{eq:limdist}) and (\ref{eq:limdistw}) of the quasi-maximum likelihood estimator and its weighted bootstrap version, also confidence sets can be constructed. For the limiting distribution of the residual bootstrap QMLE, see \citet{hy2003}. Table \ref{tab:tab2} reports the average coverage of the confidence sets, the row 'Empirical' contains the 95\% and 99\% coverage of R=1000 samples, while the other two rows show the coverage of residual and weighted bootstrap QML estimates with $R=1000$ samples and $B=1000$ bootstrap replications. Note that in each case the weighted bootstrap QMLEs performed a bit better than the residual ones. Figure \ref{fig:confset} represents the estimated pairs of parameters ($\hat{\omega},\hat{\alpha}$) and the 95\% and 99\% confidence sets -- according to the limiting distribution -- for different sample sizes (500, 1000, 2000). It can be seen that the confidence ellipses have a leaning longitudinal axis and the larger the sample size is, the smaller the ellipses become. The figures a.)--c.) were plotted 
for $R=1000$ samples and the figures d.)--f.) were plotted 
for the weighted bootstrap QMLEs, bootstrapped $B=100$ times. Compared the points against the coverage sets, the coverage looks quite decent, and there are no clusters on the outside of the ellipses.

\begin{table}[ht]
\begin{center}   
\begin{tabular}{|l|cc|cc|cc|}
\hline 
Method & \multicolumn{2}{c|}{$n=500$} & \multicolumn{2}{c|}{$n=1000$} & \multicolumn{2}{c|}{$n=1000$} \\
& 95\% & 99\% & 95\% & 99\% & 95\% & 99\% \\ \hline
 Empirical & 95.40 & 99.20 & 96.20 & 98.90 & 95.90 & 99.20 \\
 RB & 95.84 & 99.12 & 95.97 & 99.19 & 96.02 & 99.24 \\
 WB & 95.67 & 99.02 & 95.91 & 99.13 & 95.98 & 99.23 \\ \hline
\end{tabular} \end{center}
\caption{Average coverage of the 95\% and 99\% confidence sets for sample sizes 500, 1000, 2000; using residual bootstrap (RB) or weighted bootstrap (WB) methods.} 
\label{tab:tab2}
\end{table}

\section{Conclusions}

We have demonstrated that the multiplier bootstrap method reflects well the
properties of the original QMLE estimator, thus it may be used for investigating the
 estimators in practical problems (we plan to come back to this issue in another paper soon).

Another important observation of our simulations is that the asymptotic results 
presented in Sections 2 and 3 can be used for sample sizes 
in the range of thousands only, as for smaller samples the deviations may still be
substantial.

It is also worth mentioning that we have found an interesting dependence between the 
asymptotic covariance matrix and the parameter values themselves, which 
should be taken into account in practical applications.

\section{Acknowledgements}
The work was supported by the European Union Social Fund (Grant Agreement No.TÁMOP 
		4.2.1/B-09/1/KMR-2010-0003).

\section{Appendix}

\textbf{Proof of Theorem \ref{tetelboot1}.} \\
We follow the proof of \citet{fz2004} and go into details only when changes are necessary. See the original proof in their 
paper or in their book (\citet{fzbook}) on pages 156-159.

First, we introduce some notations to write the system of equations 
\begin{equation*}
\label{Bes}
\sigma_t^2 = \omega + \sum\limits_{i=1}^q \alpha_i \varepsilon^2_{t-i} + \sum\limits_{j=1}^p \beta_j \sigma^2_{t-j} \qquad 
 t \in \mathbb{Z}
\end{equation*}
in matrix form. 

\hely
$ \underline{\sigma}_t^2:=
\left( 
\begin{array}{c}
\sigma_t^2 \\ \sigma_{t-1}^2 \\ \vdots \\ \sigma_{t-p+1}^2 
\end{array} \right), \
\underline{c}_t^2:=
\left( 
\begin{array}{c}
\omega+\sum\limits_{i=1}^q \alpha_i X^2_{t-i} \\ 0 \\ \vdots \\ 0 
\end{array} \right), \
B:=
 \left( 
\begin{array}{cccc}
 \beta_1 & \beta_2 & \cdots & \beta_p \\
 1 & 0 & \cdots & 0 \\
 \vdots & \ddots & \ddots & \vdots \\
 0 & \cdots & 1 & 0
\end{array}  \right).
$

So we have 
\begin{equation} \label{sig}
\underline{\sigma}_t^2= \underline{c}_t^2 + B \underline{\sigma}_{t-1}^2 \qquad  t \in \mathbb{Z}.
\end{equation}
Let us denote by $B_k(\theta)$ the open sphere with center $\theta$ and radius $k$.

The proof consists of five steps and we also need a modification of the ergodic theorem.

\hely
\textbf{(I.) The initial values are asymptotically irrelevant} 
\begin{equation} \label{anasy1}
\underset{\theta \in \Theta}{\sup} \left| I_n^*(\theta)-\tilde{I}_n^*(\theta) \right| \xrightarrow[n \to \infty]{a.s.} 0.
\end{equation}

Iterating (\ref{sig}), we get that for some appropriate $K>0$ and $ 0< \rho <1$ 
\begin{equation} \label{iter}
\underset{\theta \in \Theta}{\sup} \norm{\underline{\sigma}_t^2(\theta)-\tilde{\underline{\sigma}}_t^2(\theta)}
\overset{a.s.}{\leq} K \rho ^t \qquad t \in \mathbb{Z}.
\end{equation}
 
\begin{flalign*}
 \underset{\theta \in \Theta}{\sup} \left| I_n^*(\theta)-\tilde{I}_n^*(\theta) \right|  \leq & 
 \frac{1}{n} \sum \limits_{t=1}^n \tau_{nt} \underset{\theta \in \Theta}{\sup} \left\{ \left| 
  \frac{\tilde{\sigma}_t^2-\sigma_t^2}{\tilde{\sigma}_t^2 \sigma_t^2} X_t^2 \right| +
  \left| \log \frac{\sigma _t^2}{\tilde{\sigma}_t^2} \right|  \right\}  \overset{(\ref{iter})}{\leq} \\ & \overset{(\ref{iter})}{\leq}
 \left(  \underset{\theta \in \Theta}{\sup} \frac{1}{\omega^2} \right) \frac{1}{n} \sum\limits_{t=1}^n \tau_{nt}
 \rho^t X_t^2 + 
 \left(  \underset{\theta \in \Theta}{\sup} \frac{1}{\omega} \right) \frac{1}{n}K \sum\limits_{t=1}^n \tau_{nt} \rho^t
\end{flalign*}
 To prove (\ref{anasy1}), it is sufficient to show that 
 \begin{equation} \label{BC1}
  \frac{1}{n} \sum\limits_{t=1}^n \tau_{nt} \rho^t X_t^2 \xrightarrow[n \to \infty]{a.s.} 0
 \end{equation} \quad and 
 \begin{equation} \label{BC2}
  \frac{1}{n} \sum\limits_{t=1}^n \tau_{nt} \rho^t \xrightarrow[n \to \infty]{a.s.} 0.
 \end{equation}
For arbitrary $\delta >0$
\begin{flalign*}
 \sum\limits_{t=0}^{\infty} P \left( \tau_{nt} \rho^t X_t^2 > \delta \right) \leq
 \sum\limits_{t=0}^{\infty} \rho^{st} \frac{E \left( \tau_{nt}^s X_t^{2s} \right)}{\delta^s} =
 \frac{E \left( \tau_{nt}^s \right)E \left(  X_t^{2s} \right) }{(1-\rho^s) \delta^s} < \infty
\end{flalign*}
and 
\begin{flalign*}
 \sum\limits_{t=0}^{\infty} P \left( \tau_{nt} \rho^t > \delta \right) \leq 
 \sum\limits_{t=0}^{\infty} \rho^t \frac{E \left( \tau_{nt}  \right)}{\delta} =
 \frac{E \left( \tau_{nt} \right) }{(1-\rho) \delta} < \infty .
\end{flalign*}
In the estimation above we applied Markov's inequality and Theorem \ref{tethat}. \\
Using the Borel-Cantelli lemma, we get
\begin{flalign*}
 P \left( \underset{t \to \infty}{\lim} \tau_{nt} \rho^t X_t^2 =0 \right) =1
\end{flalign*}
and 
\begin{flalign*}
 P \left( \underset{t \to \infty}{\lim} \tau_{nt} \rho^t =0 \right) =1.
\end{flalign*}
Finally, using Cesaro's lemma, (\ref{BC1}) and (\ref{BC2}) are proved.

\hely
\textbf{(II.) Identifiability of the parameter} 
\begin{equation*} 
 \exists t \in \mathbb{Z} \text{ such that } \quad 
 \sigma_t^2(\theta) \overset{P_{\theta_0}-a.s.}{=}  \sigma_t^2(\theta_0) \quad \Longrightarrow \quad \theta = \theta_0 .
\end{equation*}
For details, see \citet{fzbook}, page 158.

\textbf{(III.) The log likelihood function is integrable at $\theta_0$ and it has a unique minimum at the true value} 
\begin{equation*} 
 E_{\theta_0} \left| l_{nt}^*(\theta_0) \right| < \infty \text{ and if }
 \theta \neq \theta_0, E_{\theta_0} l_{nt}^*(\theta) > E_{\theta_0} l_{nt}^*(\theta_0).
\end{equation*}

It is easy to show that $E_{\theta_0}I_n^*(\theta)=E_{\theta_0}l_{nt}^*(\theta)\in \mathbb{R} \cup \{ \infty \}$,
because 
\begin{flalign*} 
 E_{\theta_0}[l_{nt}^*(\theta)]^- & = 
 E_{\theta_0} \left[ \tau_{nt} \left( \frac{X^2_t}{\sigma^2_t(\theta)}+\log \sigma^2_t(\theta) \right) \right]^- = \\ & =
 E(\tau_{nt} ) \cdot E_{\theta_0} \left( \frac{X^2_t}{\sigma^2_t(\theta)}+\log \sigma^2_t(\theta) \right)^- \leq \\ & 
 \leq E(\tau_{nt}) \cdot E_{\theta_0} \left( \log \sigma^2_t(\theta) \right)^- \leq 
  E(\tau_{nt}) \cdot E_{\theta_0} \log ^- (\omega)  < \infty.
\end{flalign*}
 
The log likelihood function is integrable at $\theta_0$: 
\begin{flalign*} 
 E_{\theta_0}l_{nt}^*(\theta_0) & = 
 E_{\theta_0} \left[ \tau_{nt} \left( \frac{\sigma^2_t(\theta_0) \eta_t^2}{\sigma^2_t(\theta_0)}+
  \log \sigma^2_t(\theta_0) \right) \right] = \\ &
 = E(\tau_{nt} ) \cdot  E_{\theta_0} \left( \eta_t^2+ \log \sigma^2_t(\theta_0)\right) =  
  E(\tau_{nt} ) \cdot \left( 1+ E_{\theta_0} \log \sigma^2_t(\theta_0) \right) < \infty
\end{flalign*}

The limit criterion is minimized at the true value  $\theta_0$ 
\begin{flalign*} 
 E_{\theta_0}l_{nt}^*(\theta) - E_{\theta_0}l_{nt}^*(\theta)_0 & \geq  
 E(\tau_{nt} ) \cdot  E_{\theta_0} \left[ \log \left( \frac{\sigma_t^2(\theta)}{\sigma_t^2(\theta_0)} \right) + 
 \log \left( \frac{\sigma_t^2(\theta_0)}{\sigma_t^2(\theta)} \right) \right] = 0.  
\end{flalign*} 
where the equality holds iff $\sigma_t^2(\theta) \overset{P_{\theta_0}-a.s.}{=} \sigma_t^2(\theta_0)$ and as a consequence of 
\textbf{(II.)}, this is equivalent to 
$\theta \overset{P_{\theta_0}-a.s.}{=} \theta_0$.
 
\textbf{(IV.) For any $\theta \neq \theta_0$, there exists a neighborhood $V(\theta)$ such that }
\begin{flalign*} 
 \underset{n \to \infty}{\text{liminf}} \underset{\breve{\theta} \in V(\theta)}{\text{inf}} \tilde{I}^*_n (\breve{\theta})
 \overset{\text{a.s.}}{>}
 E_{\theta_0}l_1(\theta_0). 
\end{flalign*}
To prove this, we use \textbf{(I.)} and a consequence of the ergodic theorem.

\begin{flalign*} 
 \underset{n \to \infty}{\text{liminf}} \underset{\breve{\theta} \in V_{1/k}(\theta) \cap \Theta}{\text{inf}} \tilde{I}^*_n (\breve{\theta})
 & \geq  
 \underset{n \to \infty}{\text{liminf}} \underset{\breve{\theta} \in V_{1/k}(\theta) \cap \Theta}{\text{inf}} I^*_n (\breve{\theta}) -
 \underset{n \to \infty}{\text{limsup}} \underset{\theta \in \Theta}{\text{ sup}} |I^*_n(\breve{\theta}) - \tilde{I}^*_n(\breve{\theta})|
 \overset{\textbf{(I.)}}{\geq} \\ & \overset{\textbf{(I.)}}{\geq} 
 \underset{n \to \infty}{\text{liminf}} \frac{1}{n} \sum\limits_{t=1}^n
  \underset{\breve{\theta} \in V_{1/k}(\theta) \cap \Theta}{\text{inf}} l^*_{nt} (\breve{\theta})=
  E_{\theta_0} \underset{\breve{\theta} \in V_{1/k}(\theta) \cap \Theta}{\text{inf}} l_1(\breve{\theta})
\end{flalign*}
In the last equation, we used  that $\underset{\theta}{\text{inf }} l^*_{nt} (\breve{\theta})$  is an ergodic process.
The expression $\underset{\breve{\theta} \in V_{1/k}(\theta) \cap \Theta}{\text{inf}} l_1(\breve{\theta})$ is monotonically increasing in $k$, so $E_{\theta_0} \underset{\breve{\theta} \in V_{1/k}(\theta) \cap \Theta}{\text{inf}} l_1(\breve{\theta})$ is also monotonically increasing and using Beppo Levi's theorem, 
\begin{flalign*} 
 E_{\theta_0} \underset{\breve{\theta} \in V_{1/k}(\theta) \cap \Theta}{\text{inf}} l_1(\breve{\theta})
 \overset{k \to \infty}{\longrightarrow} E_{\theta_0} l_1(\theta).
\end{flalign*}

\textbf{(V.) Last step of the proof, using the compactness of $\Theta$.} \\
For any neighborhood $V(\theta_0)$ of $\theta_0$, 
\begin{flalign} \label{otos}
 \underset{n \to \infty}{\text{limsup}} \underset{\breve{\theta} \in V(\theta_0)}{\text{inf}} \tilde{I}^*_n (\breve{\theta}) \leq
 \underset{n \to \infty}{\text{lim}} \tilde{I}^*_n (\theta_0) = 
 \underset{n \to \infty}{\text{lim}} I^*_n (\theta_0) = E_{\theta_0}l_1(\theta_0).
\end{flalign}
As $\Theta$ is a compact set, by definition, there exist $V(\theta_0),V(\theta_1),...,V(\theta_k)$ 
open subsets of $\mathbb{R}^{p+q+1}$, for which 
$\Theta \subseteq \left( \cup_{i=0}^k V(\theta_i) \right)$ and 
$V(\theta_1),...,V(\theta_k)$ satisfy \textbf{(IV.)}. So
\begin{flalign*}
 \underset{\theta \in \Theta}{\text{inf}} \tilde{I}^*_n (\theta) = 
 \underset{0 \leq i \leq k}{\text{min}} \ \underset{\theta \in \Theta \cap V(\theta_i)}{\text{inf}} \tilde{I}^*_n (\theta).
\end{flalign*}
As a consequence of \textbf{(IV.)} and (\ref{otos}), for $n$ large enough $\breve{\theta}_n^*$ $\in$ $V(\theta_0)$ with probability 1.
This is true for any neighborhood $V(\theta_0)$, therefore
\begin{flalign*}
\hat{\theta}^*_n \xrightarrow[n \to \infty]{a.s.} \theta_0. \ \ \Box
\end{flalign*}

\textbf{Proof of Theorem \ref{tetelboot2}.} \\
We follow the proof of \citet{fz2004} and go into details only when changes are necessary. See the original proof in their 
paper or in their book (\citet{fzbook}) on pages 159-168.\\
The Taylor-expansion of the function $\tilde{l}_{nt}^*(\theta)$ around $\theta_0$ is
\begin{flalign*} 
\tilde{l}_{nt}^*(\theta) = \tilde{l}_{nt}^*(\theta_0) + 
\frac{\partial}{\partial \theta}\tilde{l}_{nt}^*(\breve{\theta})(\theta-\theta_0),
\end{flalign*}
where $\breve{\theta}$ is between $\theta_0$ and $\theta$. \\
Derivating, summarizing and multiplying this equation with $\frac{1}{\sqrt{n}}$, we get
\begin{flalign*} 
 0  & \overset{(A5)}{=} \frac{1}{\sqrt{n}} \sum\limits_{t=1}^n \frac{\partial}{\partial \theta} \tilde{l}_{nt}^*(\hat{\theta}^*_n) 
 = \\ & \ = \frac{1}{\sqrt{n}} \sum\limits_{t=1}^n \frac{\partial}{\partial \theta} \tilde{l}_{nt}^*(\theta_0 ) + 
 \left( \frac{1}{n} \sum\limits_{t=1}^n \frac{\partial^2}{\partial \theta \partial \theta'} \tilde{l}_{nt}^*(\breve{\theta}) \right)
 \sqrt{n}(\hat{\theta}^*_n-\theta_0),
\end{flalign*}
where $\breve{\theta}$ is between $\theta_0$ and $\hat{\theta}^*_n$. \\
We will show that
\begin{flalign} 
 \label{tetel21} \frac{1}{\sqrt{n}} \sum\limits_{t=1}^n \frac{\partial}{\partial \theta} \tilde{l}_{nt}^*(\theta_0 )  
 \xrightarrow[n \to \infty]{d}
 N \left( 0,\gamma(\kappa_{\eta}-1)J \right) \\
 \label{tetel22} \frac{1}{n} \sum\limits_{t=1}^n \frac{\partial^2}{\partial \theta_i \partial \theta_j} \tilde{l}_{nt}^*(\breve{\theta})
 \xrightarrow[n \to \infty]{d} J(i,j).
\end{flalign} 
 
The proof consists of six steps.

\hely
\textbf{(I.) Integrability of the second-order derivatives of $l_{nt}^*(\theta)$ at $\theta_0$} 
\begin{flalign*}
 E_{\theta_0} \left \| \frac{\partial^2 l_{nt}^*}{\partial \theta \partial \theta^T}(\theta_0)  \right \| < \infty .
\end{flalign*}
As $E\tau_{nt} < \infty$ and $\tau_{nt}$ is independent from $l_t(\theta_0)$, it is sufficient to show that
\begin{flalign*}
 E_{\theta_0} \left \| \frac{\partial^2 l_t}{\partial \theta \partial \theta^T}(\theta_0)  \right \| < \infty ,
\end{flalign*}
which is proven in \citet{fzbook}, on pages 160-162.

\hely
\textbf{(II.) J is invertible and $\text{Var}_{\theta_0} \left( \frac{\partial l_{nt}^*}{\partial \theta}(\theta_0) \right) = 
E \tau^2_{nt} \cdot (\kappa_{\eta}-1)J$} \\
The invertibility of $J$ is verified in \citet{fzbook}, on page 163. \\
Using (I.), $E\tau_{nt} <\infty$ and the independence between $\tau_{nt}$ and $l_t(\theta)$, we have
\begin{flalign*}
 E_{\theta_0} \left( \frac{\partial l_{nt}^*}{\partial \theta}(\theta_0) \right) = 
  E \tau_{nt} \cdot \underbrace{E_{\theta_0}(1-\eta_t^2)}_{1-1=0} \cdot
  E_{\theta_0} \left( \frac{1}{\sigma_t^2(\theta_0)} \cdot \frac{\partial \sigma_t^2}{\partial \theta} (\theta_0) \right) =0
\end{flalign*}
Then we obtain
\begin{flalign*}
 \text{Var}_{\theta_0} \left( \frac{\partial l_{nt}^*}{\partial \theta}(\theta_0) \right) & = 
 E_{\theta_0} \left( \frac{\partial l_{nt}^*}{\partial \theta}(\theta_0) \cdot
   \frac{\partial l_{nt}^*}{\partial \theta ^T}(\theta_0) \right) = \\
 & = E \tau_{nt}^2 \cdot \underbrace{E_{\theta_0}(1-\eta_t)^2}_{\kappa_{\eta}-1} \cdot
  E_{\theta_0} \left( \frac{1}{\sigma_t^4(\theta_0)} \cdot \frac{\partial \sigma_t^2}{\partial \theta} (\theta_0) 
  \cdot \frac{\partial \sigma_t^2}{\partial \theta^T} (\theta_0)  \right) = \\
 & = E \tau_{nt}^2 \cdot (\kappa_{\eta}-1) \cdot J.
\end{flalign*}

\hely
\textbf{(III.) Uniform integrability of the third-order derivatives of $l_{nt}^*(\theta)$ at $\theta_0$:} \\
There exists a neighborhood $V(\theta_0)$ of $\theta_0$ such that, for all $i,j,k \in \{ 1,...,p+q+1 \}$,
\begin{flalign*}
 E_{\theta_0} \underset{\theta \in V(\theta_0)}{\sup} 
  \left| \frac{\partial^3l_{nt}^*(\theta)}{\partial \theta_i \partial \theta_j \partial \theta_k} \right| < \infty.
 \end{flalign*}
As $E\tau_{nt} < \infty$ and $\tau_{nt}$ is independent from $l_t(\theta_0)$, it is sufficient to show that
\begin{flalign*}
 E_{\theta_0} \underset{\theta \in V(\theta_0)}{\sup} 
  \left| \frac{\partial^3l_t(\theta)}{\partial \theta_i \partial \theta_j \partial \theta_k} \right| < \infty,
 \end{flalign*}
which is proven in \citet{fzbook}, on pages 163-165.

\hely
\textbf{(IV.) The initial values are asymptotically irrelevant:} 
\begin{flalign}
 \label{IV1} \left \| \frac{1}{\sqrt{n}} \sum\limits_{t=1}^n \left(  \frac{\partial l_{nt}^*}{\partial \theta}(\theta_0) -
  \frac{\partial \tilde{l}_{nt}^*}{\partial \theta}(\theta_0) \right) \right \| \xrightarrow[n \to \infty]{p} 0 \quad \text{and} \\
 \label{IV2}  \underset{\theta \in V(\theta_0)}{\sup} 
  \left \| \frac{1}{n} \sum\limits_{t=1}^n \left(  \frac{\partial^2 l_{nt}^*}{\partial \theta \partial \theta^T}(\theta) -
  \frac{\partial^2 \tilde{l}_{nt}^*}{\partial \theta \partial \theta^T}(\theta) \right) \right \| \xrightarrow[n \to \infty]{p} 0.
\end{flalign}
Using the results of \citet{fzbook} (pages 165-166) we have
\begin{flalign*}
\left| \frac{\partial l_{nt}^*}{\partial \theta_i}(\theta_0) - 
  \frac{\partial \tilde{l}_{nt}^*}{\partial \theta_i}(\theta_0) \right| \leq
  K \tau_{nt} \rho^t(1+\eta _t^2) \left| 1+ 
   \frac{1}{\sigma_t^2(\theta_0)} \cdot \frac{\partial \sigma_t^2}{\partial \theta_i} (\theta_0) \right|.
\end{flalign*}
So we obtain the estimate
\begin{flalign*}
 \frac{1}{\sqrt{n}} \sum\limits_{t=1}^n 
 \left|  \frac{\partial l_{nt}^*}{\partial \theta_i}(\theta_0) - \frac{\partial \tilde{l}_{nt}^*}{\partial \theta_i}(\theta_0) \right|
 \leq \breve{K} \frac{1}{\sqrt{n}} \sum\limits_{t=1}^n  \tau_{nt} \rho^t(1+\eta _t^2) \left| 1+ 
  \frac{1}{\sigma_t^2(\theta_0)} \cdot \frac{\partial \sigma_t^2}{\partial \theta_i} (\theta_0) \right|.
\end{flalign*}
Markov's inequality, the independence between $\tau_{nt}$, $\eta_t$ and $\sigma_t^2(\theta_0)$ imply that, for
all $\varepsilon >0$,
\begin{flalign*}
 P \left( \frac{1}{\sqrt{n}} \sum\limits_{t=1}^n  \tau_{nt} \rho^t(1+\eta _t^2) \left| 1+ 
  \frac{1}{\sigma_t^2(\theta_0)} \cdot \frac{\partial \sigma_t^2}{\partial \theta_i} (\theta_0) \right| > \varepsilon \right) \leq \\
 \leq \frac{2}{\sqrt{n} \varepsilon} \left( 1+ E_{\theta_0}\left|  
  \frac{1}{\sigma_t^2(\theta_0)} \cdot \frac{\partial \sigma_t^2}{\partial \theta_i} (\theta_0) \right| \right)
  \sum\limits_{t=1}^n  \rho^t E \tau_{nt}  ,
\end{flalign*}
where $0<\rho<1$. \\
To show (\ref{IV1}), it is sufficient to prove that $\underset{n \to \infty}{\lim} \sum\limits_{t=1}^{n}  
 \rho^t E \tau_{nt} < \infty$: 
\begin{flalign*}
 \underset{n \to \infty}{\lim} \sum\limits_{t=1}^{n}   \rho^t E \tau_{nt} \overset{\textbf{B3}}{=}
 \underset{n \to \infty}{\lim}  E \tau_{n1} \sum\limits_{t=1}^{n}   \rho^t = \frac{\rho}{1-\rho} < \infty.
\end{flalign*}
(\ref{IV2}) can be proven similarly.

\hely
\textbf{(V.) } Using the martingale CLT (or Lindeberg's CLT), we prove that
\begin{flalign} \label{V}
\frac{1}{\sqrt{n}} \sum\limits_{t=1}^n  \frac{\partial l_{nt}^*}{\partial \theta}(\theta_0)  \xrightarrow[n \to \infty]{d} 
 N \left(0 , \gamma  (\kappa_{\eta}-1)J) \right).
\end{flalign}
Let $\mathcal{F}_{nt}=\mathcal{F}_{t}=\sigma ( \{ X_t,X_{t-1},... \} )$ and \\ 
for all $\lambda \in \mathbb{R}^{p+q+1}$ \ 
 $\eta_{nt}=\frac{1}{\sqrt{n}} \lambda^T\frac{\partial l_{nt}^*}{\partial \theta}(\theta_0) = 
 \frac{\tau_{nt}}{\sqrt{n}} \lambda^T\frac{\partial l_t}{\partial \theta}(\theta_0) $.  \\
So for every $n$, $(\eta_{nt},\mathcal{F}_{nt})_{t\in \mathbb{Z}}$ is a square integrable martingale difference. \\
Let us denote with $\sigma^2_{nt}=E_{\theta_0}(\eta_{nt}^2|\mathcal{F}_{t-1})$, therefore the process 
\begin{flalign*} 
(\sigma^2_{nt})_{t=1,...,n}= \frac{1}{n} 
  \left[ E_{\theta_0}\left( \left. \tau_{nt}^2 \left[ \lambda^T\frac{\partial l_t}{\partial \theta}(\theta_0)
\right]^2 \right| \mathcal{F}_{t-1} \right) \right]_{t=1,...,n} 
\end{flalign*}
is stationary and ergodic. \\
As a consequence, using \textbf{B6} for Bernstein's theorem
\begin{flalign*} 
\sum\limits_{t=1}^{n}\sigma^2_{nt} & =\frac{1}{n} \sum\limits_{t=1}^{n} 
E_{\theta_0}\left( \left. \tau_{nt}^2 \left[ \lambda^T \frac{\partial l_t}{\partial \theta}(\theta_0) \right]^2
\right| \mathcal{F}_{t-1} \right) \xrightarrow[n \to \infty]{p} \\ & \xrightarrow[n \to \infty]{p}
E_{\theta_0} \left[ E_{\theta_0}\left( \left. \underset{n \to \infty}{\lim} \tau_{n1}^2 \left[ \lambda^T\frac{\partial l_1}{\partial \theta}(\theta_0) \right]^2
\right| \mathcal{F}_{t-1} \right) \right] = \gamma \cdot (\kappa_{\eta}-1) \cdot J.
\end{flalign*} 
We also have for all $\varepsilon>0$
\begin{flalign*} 
\sum\limits_{t=1}^{n} E_{\theta_0}  \left[ \eta^2_{nt}I(|\eta_{nt}| \geq \varepsilon) \right]=
\sum\limits_{t=1}^{n} \frac{1}{n} \int_{ \left\{ \left| \tau_{nt} \lambda^T\frac{\partial l_t}{\partial \theta}(\theta_0)
 \right| \geq \sqrt{n} \varepsilon \right\}} 
\tau_{nt}^2 \left[ \lambda^T\frac{\partial l_t}{\partial \theta}(\theta_0) \right]^2 dP_{\theta_0} = \\
= \int_{ \left\{ \left| \tau_{n1} \lambda^T\frac{\partial l_1}{\partial \theta}(\theta_0)
 \right| \geq \sqrt{n} \varepsilon \right\}} 
\tau_{n1}^2 \left[ \lambda^T\frac{\partial l_1}{\partial \theta}(\theta_0) \right]^2 dP_{\theta_0}
\xrightarrow{n \to \infty} 0.
\end{flalign*}
At the second equality we used the stationarity of the process. \\
Using the martingale CLT on the process $(\eta_{nt},\mathcal{F}_{nt})_{t\in \mathbb{Z}}$ and then the
Cramér-Wold theorem, (\ref{V}) is proved.
 
\hely
\textbf{(VI.)} Using the second order derivative of the Taylor expansion of $l_{nt}^*$, it can be seen that 
\begin{flalign*}
\frac{1}{n} \sum\limits_{t=1}^n  \frac{\partial l_{nt}^*}{\partial \theta_i \partial \theta_j}(\breve{\theta}_{ij})  \xrightarrow[n \to \infty]{a.s.}  J(i,j).
\end{flalign*}
At last, if we combine \textbf{(IV.)}, \textbf{(V.)},\textbf{(VI.)} and apply Slutsky's lemma on the first order derivative of 
the Taylor expansion of $l_{nt}^*$, (\ref{tetelboot2}) is proved. \ \ $\Box$

\bibliographystyle{plainnat}
\bibliography{refs}

\end{document}